\documentclass[12pt]{article}
\usepackage{amsmath,amsthm,amsfonts,latexsym,epsfig}
\newtheorem{theorem}{Theorem}[section]
\newtheorem{lemma}[theorem]{Lemma}
\newtheorem{proposition}[theorem]{Proposition}
\title{Tree-like properties of cycle factorizations}
\author{Ian Goulden\footnote{Department of Combinatorics and
Optimization, University of Waterloo, email: {\small\texttt
ipgoulden@math.uwaterloo.ca}}~
and Alexander Yong\footnote{Department of Mathematics,
University of Michigan, Ann Arbor, email: {\small\texttt
ayong@umich.edu}}}
\date{\today}
\begin{document}
\maketitle

\begin{abstract}
\noindent
We provide a bijection between the set of {\em factorizations}, that is, ordered $(n-1)$-tuples
of transpositions in ${\mathcal S}_{n}$ whose product is $(12...n)$, and labelled trees on $n$
vertices. We prove a refinement of a theorem of D\'{e}nes [\ref{D}] that
establishes new tree-like
properties of factorizations. In particular, we show that a certain class
of transpositions of a
factorization correspond naturally under our bijection to leaf edges of a tree. Moreover, we
give a generalization of this fact.
\end{abstract}

\section{Introduction}

	Let ${\mathcal T}_{n}$ be the set of labelled trees on vertices
$\{1,2,...,n\}=[n]$, and ${\mathcal F}_{n}$ be the set of $(n-1)$-tuples
of transpositions $(\sigma_{1},\ldots ,\sigma_{n-1})$ in the symmetric
group ${\mathcal S}_{n}$ acting on $[n]$, whose ordered product
$\sigma_{1}\cdots\sigma_{n-1}$ is equal to the cycle ${\mathcal C}_{n}=(12\ldots n)$.
The elements of ${\mathcal F}_{n}$ are called {\em factorizations}, and the
transpositions $\sigma_{1},\ldots ,\sigma_{n-1}$ in a factorization are called
{\em factors}. Cayley [\ref{C}] proved that $\mid {\mathcal T}_{n}\mid=n^{n-2}$,
and D\'{e}nes [\ref{D}] proved that $\mid {\mathcal F}_{n}\mid=\mid {\mathcal T}_{n}\mid$,
by giving a bijection between sets of cardinality $(n-1)!\mid {\mathcal F}_{n}\mid$
and $(n-1)!\mid {\mathcal T}_{n}\mid$. D\'{e}nes posed the problem of
finding a bijection between ${\mathcal F}_{n}$ and ${\mathcal T}_{n}$, and
subsequently two such bijections have been given, by  Moszkowski [\ref{M}]
and Goulden and Pepper [\ref{GP}].

\bigskip
	Although both of these bijections are reasonably simple, neither of them
restricts nicely to natural combinatorial subsets (e.g. so that the image of
a combinatorially natural subset of ${\mathcal T}_{n}$ corresponds to a
natural subset of ${\mathcal F}_{n}$). However, by examining the elements
of
${\mathcal T}_{n}$ and ${\mathcal F}_{n}$ for small $n$, we find that there
are natural combinatorial subsets of ${\mathcal T}_{n}$ and ${\mathcal F}_{n}$
of equal cardinality, as follows. Let ${\mathcal T}_{n}(k)$ be the set of
trees in ${\mathcal T}_{n}$ with $k$ leaves (vertices of degree one). A
transposition $(s \ t)$ on $[n]$ is called a {\em consecutive pair} if
$t \equiv s+1$ modulo $n$, where throughout, we write $n$ to mean 0. Let ${\mathcal
F}_{n}(k)$
be the set of factorizations in ${\mathcal F}_{n}$ with $k$ factors that are
consecutive pairs. Table 1 gives the cardinalities $\mid{\mathcal T}_{n}(k)\mid$
for $n\leq 6$, and a systematic examination of the factorizations in ${\mathcal F}_{n}$ shows that
$\mid {\mathcal F}_{n}(k)\mid=\mid {\mathcal T}_{n}(k)\mid$ for $3\leq
n\leq 6$. This suggests
that
there exists a bijection between ${\mathcal F}_{n}$ and ${\mathcal T}_{n}$
for arbitrary $n\geq 3$ that maps consecutive pairs to leaves, but neither
of the previous bijections exhibits this property.

\begin{center}
\begin{tabular}{|l||l|l|l|l|l|} \hline
$n/k$ & 2 & 3 & 4 & 5 & total \\ \hline\hline
2   & 1   &     &     &   & 1   \\ \hline
3   & 3   &     &     &   & 3    \\   \hline
4   & 12  & 4   &     &   & 16   \\   \hline
5   & 60  & 60  & 5   &   & 125  \\   \hline
6   & 360 & 720 & 210 & 6 & 1296 \\   \hline
\end{tabular}
\end{center}
\begin{center}
\bigskip
Table 1: The number of trees on $n$ vertices
with $k$ leaves for $n\leq 6$.
\end{center}

\bigskip
	In this paper we describe a bijection between ${\mathcal F}_{n}$
and ${\mathcal T}_{n}$ in which consecutive pairs of factorizations correspond
to leaves of trees. We refer to our bijection as a {\em structural} bijection
because of this correspondence between these combinatorial structures
(consecutive pairs and leaves). But more is true; the bijection extends
to generalizations of consecutive pairs and of leaves respectively, as
described below.

\bigskip
	For a tree $T\in {\mathcal T}_{n}$, consider removing any single
edge from the tree, to get two trees $T_{1}$ and $T_{2}$ (the
components of the graph that results when the edge is deleted from $T$).
Let $t_{i}$, $i=1,2$, be the number of vertices
in $T_{i}$ (so, e.g., $t_{1}+t_{2}=n$),
and define the {\em edge-deletion index} of the edge to be $\min\{t_{1},t_{2}\}$.
Define the {\em edge-deletion distribution} of the tree $T$ to be
$\varepsilon(T)=(a_{1},a_{2},...)$ where $a_{j}$ is the number of edges in $T$
with edge-deletion index $j$ (so e.g., $a_{1}+a_{2}+...=n-1$). Let
${\mathcal T}_{n}(a_{1},a_{2},...)$ be the set of trees in ${\mathcal T}_{n}$
with edge-deletion distribution $(a_{1},a_{2},...)$.

\bigskip
	For a transposition $(s \ t)$, $s<t$, define the {\em difference index} to
be $\min\{t-s,n-t+s\}$. For a factorization $F\in {\mathcal F}_{n}$,
define the {\em difference distribution} of $F$ to be $\delta(F)=(d_{1},d_{2},...)$
where $d_{j}$ is the number of factors in $F$ with difference index $j$
(so e.g., $d_{1}+d_{2}+...=n-1$). Let ${\mathcal F}_{n}(d_{1},d_{2},...)$
be the set of trees in ${\mathcal F}_{n}$ with difference distribution
$(d_{1},d_{2},...)$.

\bigskip
	Our structural bijection, described in section 3, actually gives a
bijection between ${\mathcal F}_{n}(c_{1},c_{2},...)$ and
 ${\mathcal T}_{n}(c_{1},c_{2},...)$
for all $(c_{1},c_{2},...)$, $c_{i}\geq 0$, $c_{1}+c_{2}+...=n-1$, $n\geq 1$.
The main result of our paper is as follows:

\begin{theorem}
\label{main_theorem}
For each $n\geq 1$, there is a bijection

$$\phi:{\mathcal F}_{n}\to {\mathcal T}_{n}: F\mapsto T$$

\noindent
such that $\delta(F)=\varepsilon(T)$.
\end{theorem}

	Note that for factorizations, a factor with difference
index $1$ is a consecutive
pair, and for trees, an edge with edge-deletion index $1$ is
incident with a leaf,
so this is a generalization of the consecutive pair-leaf correspondence, as
promised. In the case $n=2$, the single edge in the unique tree has difference index
equal to $1$, but is incident with two leaves. Of course, the single factor
$(1 \ 2)$ in the unique factorization is a consecutive pair, so the
bijection claimed in Theorem~\ref{main_theorem} holds for $n=2$, where
the consecutive pair-leaf correspondence breaks down.

\bigskip
	The bijection is based on a geometrical interpretation of a factorization,
called a {\em chord diagram}, whose properties are developed in section 2.
The bijection, described in section 3, has a smooth composition
with the well-known Pr\"{u}fer code
bijection between trees and the set $[n]^{n-2}$. Consequently one can
obtain a bijection under this composition that canonically proves that
$\mid {\mathcal F}_{n}\mid = n^{n-2}$.

\bigskip
	Further motivation for this paper, which gives a third bijection for
the D\'{e}nes result, beyond the combinatorial benefits of exhibiting tree
properties of edges as differences of factors, is provided in section 4,
where we describe recent work on more general factorization questions in the
symmetric group, related to certain problems arising from algebraic geometry.

\bigskip
	Finally, there are immediate enumerative consequences of our main
result. For example, there is a nice formula for the
entries in Table 1, which can
be obtained in various ways by counting trees with a given number of leaves.
This formula is a simple multiple of a Stirling number of the second kind,
and in closed form it gives
$$\mid {\mathcal T}_{n}(k)\mid = {n \choose k} \sum_{i=0}^{n-k}
{n-k \choose i} (-1)^{n-k-i} i^{n-2},$$
for $2\leq k\leq n$ (see, e.g., Stanton and White~[\ref{SW}], p. 67).
Using our bijection, we therefore have established that this formula
also holds for $\mid {\mathcal F}_{n}(k)\mid$.

\section{The circle chord diagram}

	We begin with a detailed analysis of several aspects of factorizations
in ${\mathcal F}_{n}$. First, if $\rho\in {\mathcal S}_{n}$ and $\tau=(s \ t)\in {\mathcal S}_{n}$
is a transposition, then there are two cases that arise in determining
the product $\alpha=\tau\rho$. If $s,t$ appear on the same cycle in the
disjoint cycle representation of $\rho$, then that cycle is cut into
two different cycles (one containing $s$, the other $t$) in the
disjoint cycle representation of $\alpha$. Otherwise, if $s,t$ appear
on two different cycles of $\rho$, then these cycles are joined into one
cycle (containing both $s$ and $t$) of $\alpha$. We call $\tau$ a {\em join}
or a {\em cut} of $\rho$, respectively, in these cases.

\bigskip
	For a factorization $F=(\sigma_{1},...\sigma_{n-1}$) in ${\mathcal F}_{n}$,
let $f_{i}=\sigma_{i}\sigma_{i+1}...\sigma_{n-1}$, $i=1,...,n-1$, be the
{\em partial products} of $F$, and define $\sigma_{i}$ to be a {\em join} or
{\em cut} of $F$ when $\sigma_{i}$ is a join or cut of $f_{i+1}$,
respectively (when determining the product $f_{i}=\sigma_{i}f_{i+1}$),
for $i=1,...,n-1$ (where $f_{n}=e$, the identity of ${\mathcal S}_{n}$).
Now $f_{n}$ has $n$ cycles (all fixed points), and
$f_{1}={\mathcal C}_{n}$ has 1 cycle. Moreover, each join decreases the number
of cycles by 1, and each cut increases the number of cycles by 1. We conclude
that each of the $n-1$ factors $\sigma_{i}$ in
$f_{1}=\sigma_{1}...\sigma_{n-1}f_{n}$
must be a join, since together they decrease the $n$ cycles of $f_{n}$ by $n-1$,
to the single cycle of $f_{1}$.

\bigskip
	We say that a sequence $\alpha_{1},...,\alpha_{m}$ of elements in
$[n]$ is ${\mathcal C}_{n}${\em -ordered} if the order of the elements
is
consistent with their circular order on the cycle ${\mathcal C}_{n}=(12...n)$.
Equivalently, this means that there is a
unique $i$ with $1\leq\alpha_{i}<\alpha_{i+1}<...<\alpha_{m}<\alpha_{1}<...<\alpha_{i-1}\leq n.$

\begin{proposition}
\label{prop_one}
For $F\in {\mathcal F}_{n}$, and a partial product $f_{i}$ of $F$, any
subsequence of elements on a cycle of $f_{i}$ is ${\mathcal C}_{n}$-ordered.
\end{proposition}

{\noindent \bf Proof:} For $F=(\sigma_{1},...,\sigma_{n-1})$, we have
$f_{i}=\sigma_{i}...\sigma_{n-1}$, where all factors $\sigma_{1},...,\sigma_{n-1}$
are joins, from the above join-cut analysis. But $f_{1}=\sigma_{1}...\sigma_{n-1}={\mathcal C}_{n}$,
and the effect of the sequence of subsequent joins $\sigma_{1},...,\sigma_{i-1}$
in $f_{1}=\sigma_{1}...\sigma_{i-1}f_{i}$, on the elements of a cycle of
$f_{i}$ is to keep them together on cycles that are formed by the joins,
and to maintain their circular order around such cycles. We conclude that
the elements on each cycle of $f_{i}$ must be ${\mathcal C}_{n}$-ordered,
and therefore so must all subsequences of elements on each such cycle.\qed

\bigskip
	We now consider a {\em circle chord diagram}. For any fixed $n$,
this is a circle drawn in the plane with $n$ points on it, labelled
$1,2,...,n$ clockwise. In addition, there are $n-1$ chords on these $n$
points,
numbered $2,...,n$ distinctly. For example, Figure 1 gives a circle chord diagram
with $n=9$; the numbers on the edges are circled to distinguish them from
the names of points on the circle.

\begin{figure}[h]
\centering
\epsfig{file=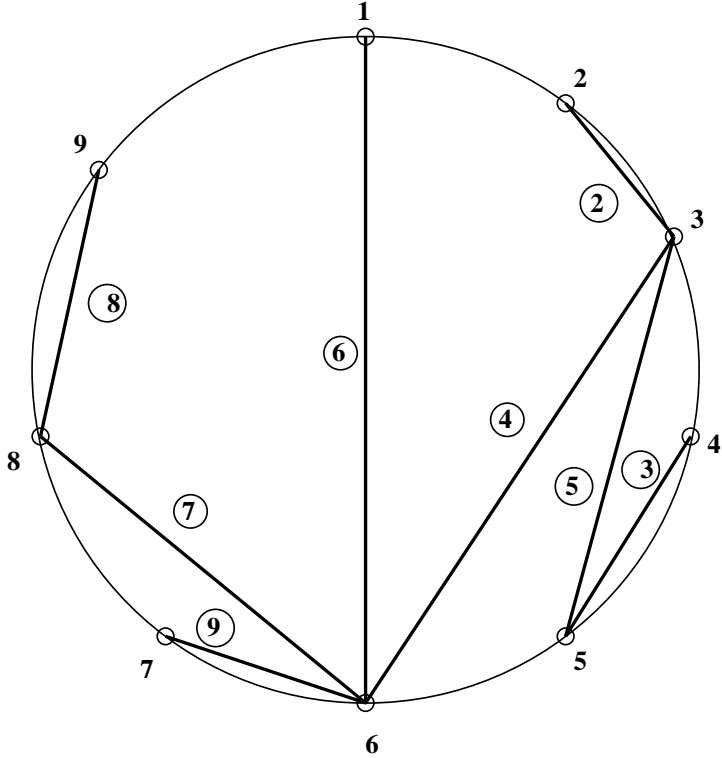,height=10cm}
\caption{A circle chord diagram with $n=9$.}
\label{figure_1}
\end{figure}

\bigskip
	There is a natural injection from factorizations to circle
chord diagrams: for $F=(\sigma_{1},...,\sigma_{n-1})\in {\mathcal F}_{n}$,
the factor $\sigma_{i}=(s_{i} \ t_{i})$ corresponds to a chord
numbered $i+1$, joining points $s_{i}$ and $t_{i}$, for $i=1,...,n-1$. Let
$C(F)$ be the circle chord diagram associated with $F$ in this way. For
example, the chord diagram illustrated in Figure 1 is $C(F_{0})$, where
\begin{equation}\label{factex}
F_{0}=((2 \ 3),(4 \ 5),(3 \ 6),(3 \ 5),(1 \ 6),(6 \ 8),(8 \ 9),(6 \ 7)).
\end{equation}

	The circle chord diagram associated with a factorization $F$ satisfies a
number of conditions, and we now establish some of these.

\begin{theorem}
\label{theorem_two}
In the circle chord diagram $C(F)$ of a factorization $F\in{\mathcal F}_{n}$,\\

\noindent
(i) the chords form a tree on $[n]$, $(*_{1})$\\
\noindent
(ii) the chords meet only at endpoints, $(*_{2})$\\
\noindent
(iii) the edge numbers on the chords encountered when
moving around a vertex clockwise
across the interior of the circle, form a decreasing sequence of elements
in $\{2,...,n\}$. $(*_{3})$
\end{theorem}

{\noindent \bf Proof:} For (i), let ${\mathcal G}_{j}$ be the graph on
vertex-set $[n]$, whose edges are the chords corresponding to factors
$\sigma_{j},...,\sigma_{n-1}$, for $j=1,...,n$ (${\mathcal G}_{n}$ has no
edges). Then ${\mathcal G}_{n}$ has $n$ components (each a single vertex),
and the condition, established above, that $\sigma_{j}$ is a join for
each $j$ implies that the chord corresponding to $\sigma_{j}$ is incident
with vertices in different components of ${\mathcal G}_{j+1}$, for each
$j=1,...,n-1$. Thus ${\mathcal G}_{j}$ has one fewer components than
${\mathcal G}_{j+1}$ for each $j=1,...,n-1$, and we conclude that ${\mathcal G}_{1}$
has one component, so it is a connected graph. But the edges of $G_{1}$
are the chords of $C(F)$, so the $n-1$ chords of $C(F)$ are a connected graph on
$n$ vertices, which must therefore be a tree.

\bigskip
	For (ii), suppose otherwise, that the chords corresponding to
$\sigma_{i}=(s \ t)$ and $\sigma_{j}=(u \ v)$, where $s<t,\, u<v$ and $i<j$,
cross each
other. Now the geometric crossing condition is equivalent to the
condition that the sequence $stuv$ is not ${\mathcal C}_{n}$-ordered.
But in $f_{i}$, the cycle containing $s$ will include $stuv$ as a
subsequence,
and we have a contradiction of Proposition \ref{prop_one}, which establishes
that $stuv$ must be ${\mathcal C}_{n}$-ordered. We conclude that chords
do not cross and can therefore meet only at endpoints.

\bigskip
	For (iii), for each fixed $i=1,...,n$, suppose the factors moving
$i$ are

$$\sigma_{l_{1}}=(i \ s_{1}),...,\sigma_{l_{k}}=(i \ s_{k})$$

\bigskip
\noindent
where $1\leq l_{1}<...<l_{k}\leq n-1$, $k\geq 1$.

\bigskip
	Then $f_{l_{1}}$ will include $is_{k}...s_{1}$ as a subsequence
on the cycle containing $i$, and we conclude from Proposition \ref{prop_one}
that $is_{k}...s_{1}$ is ${\mathcal C}_{n}$-ordered. But the
edge corresponding to $\sigma_{l_{j}}$ has number $l_{j}+1$, and (iii)
follows. \qed

\bigskip
	For example, it is straightforward to verify that the circle
chord diagram $C(F_{0})$ illustrated in Figure 1 does indeed satisfy
conditions $(*_{1}),(*_{2})$ and $(*_{3})$.

\bigskip
	Now for circle chord
diagrams satisfying condition $(*_{2})$, the $n-1$ chords and the circle
partition the circle and its interior into $n$ {\em regions}. The {\em boundary}
of the region consists of a collection of chords and {\em arcs} of the circle.
An arc is a segment of the circle from point $i$ to point $i+1$ modulo $n$.

\begin{proposition}
For circle chord diagrams satisfying conditions $(*_{1})$ and $(*_{2})$,
each region contains precisely one arc in its boundary.
\end{proposition}

{\noindent \bf Proof:} If the boundary of any region consists entirely
of chords, then these chords form a cycle in the graph of the chords (called
${\mathcal G}_{1}$ in the proof of Theorem \ref{theorem_two}(i)). But this
graph
is a tree, from Theorem \ref{theorem_two}(i), and therefore has no cycles.
We
conclude that the boundary of each of the $n$ regions contains at least
one of the $n$ arcs. But this means that each region has exactly one arc,
giving the result.\qed

\bigskip
	For example, each region of the circle chord
diagram $C(F_{0})$ illustrated in Figure 1 contains precisely one arc
in its boundary.

\bigskip
	Now consider the following condition for the above regions:
the numbers on the chords of the boundary increase clockwise, starting immediately
after the unique arc. $(*_{3})'$

\bigskip
	Note, for example, that each of the 9 regions in Figure 1
satisfies $(*_{3})'$.

\begin{proposition}
For circle chord diagrams, conditions $(*_{1}),(*_{2})$ and $(*_{3})$
are equivalent to $(*_{1}),(*_{2})$ and $(*_{3})'$.
\end{proposition}

{\noindent \bf Proof:} Immediate.\qed

\bigskip
	We end this section by showing that
conditions $(*_{1}),(*_{2})$ and $(*_{3})'$
characterize circle chord diagrams associated with factorizations.

\begin{lemma}
A circle chord diagram on $n$ points satisfying conditions $(*_{1}),(*_{2})$
and $(*_{3})'$ is equal to $C(F)$ for some $F\in {\mathcal F}_{n}$.
\end{lemma}

{\noindent \bf Proof:} Consider a circle chord diagram satisfying conditions
$(*_{1}),(*_{2})$ and $(*_{3})'$. Suppose that the chord numbered $i$
joins
points $a_{i}$ and $b_{i}$, and let $\sigma_{i-1}=(a_{i} \ b_{i})$, for
$i=2,...,n$. Now consider the product of transpositions

$$\sigma=\sigma_{1}...\sigma_{n-1}.$$

\bigskip
\noindent
Condition $(*_{3})'$ implies that $\sigma(j)\equiv j+1$ modulo $n$ for
each $j=1,...,n$, by considering the action of the transpositions on the
boundary of the region containing the arc $(j, j+1)$. Thus $\sigma={\mathcal C}_{n}$,
and $F'=(\sigma_{1},...,\sigma_{n-1})$ is a factorization in ${\mathcal F}_{n}$.
The result follows, since we have established that the circle chord
diagram is equal to $C(F')$.\qed

\bigskip
	In summary, the results in this section have established that there
is a bijection between ${\mathcal F}_{n}$ and circle chord
diagrams satisfying $(*_{1}),(*_{2})$
and $(*_{3})'$. 

\section{The structural bijection}

	We are now able to describe the structural
bijection that proves our main theorem.
Consider the circle chord diagram $C(F)$ for some factorization
$F\in {\mathcal F}_{n}$. Form the graph $\phi(F)$ by a ``planar dual''
construction, as follows. For each
region of $C(F)$ we have a vertex of $\phi(F)$ (say, drawn in the ``middle''
of the region). Then place an edge between two vertices if the
boundaries of their
corresponding regions share a chord. We (temporarily) assign label $i$
to this edge of $\phi(F)$, where $i$ is the number of the shared chord in
$C(F)$. Thus, at this stage, $\phi(F)$ has $n$ vertices and $n-1$ edges
(one edge for each edge of $C(F)$), and is connected because $C(F)$ is
connected, so we conclude that $\phi(F)$ is a tree. Note that
condition $(*_{3})'$ on $C(F)$ implies that, at each vertex
of $\phi(F)$, the clockwise sequence of labels on the incident edges
is ${\mathcal C}_{n}$-ordered.

\bigskip
	Now complete the construction by labelling the vertices distinctly
with the elements of $[n]$, as follows. The vertex corresponding to the
region with arc $(n,1)$ in its boundary has label 1. For each edge,
find the unique path to vertex 1 from that edge, and ``slide'' the
temporary label on the edge to the incident vertex away from vertex 1,
thus labelling the other $n-1$ vertices $2,...,n$. The resulting
tree is $\phi(F)$, and it is clear from our
description above that $\phi(F)\in {\mathcal T}_{n}$. For example, Figure
2 illustrates $\phi(F_{0})$ where $F_0$ is
given in~(\ref{factex}), and $C(F_{0})$
is given in Figure 1.

\begin{figure}[h]
\centering
\epsfig{file=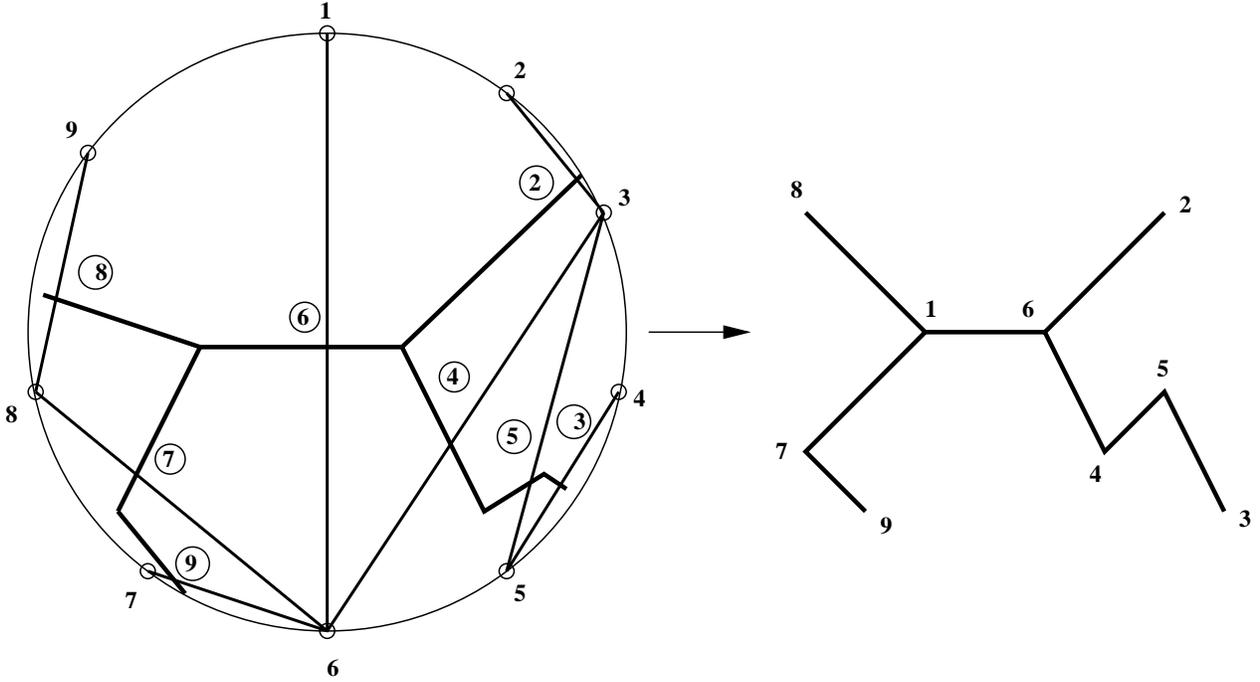,height=9cm}
\caption{Construction of the tree $\phi(F_{0})$.}
\label{figure_2}
\end{figure}

\bigskip
	We claim that $\phi:{\mathcal F}_{n}\to {\mathcal T}_{n}$ is
a bijection. In section 2 we proved that there is a bijection between
${\mathcal F}_{n}$ and circle chord diagrams satisfying $(*_{1}),(*_{2})$
and $(*_{3})'$. Also, our ``planar dual'' construction above is
a bijection between circle chord diagrams
satisfying $(*_{1}),(*_{2})$ and $(*_{3})'$ and ${\mathcal T}_{n}$,
since the ${\mathcal C}_{n}$-ordered requirement at each vertex
forces a unique planar embedding of a tree. Together, these bijections
prove the claim. In the resulting bijection $\phi$, note that
a factor in $F$ with difference index $k$
corresponds precisely to an edge of $\phi(F)$ with edge-deletion
index $k$, so $\delta(F)=\varepsilon(\phi(F))$ and we have proved Theorem
\ref{main_theorem}.

\bigskip
	To reverse the bijection, consider an arbitrary tree $T$. Now
``slide'' the label on each vertex $2,...,n$ to the incident edge along the
unique path to vertex $1$. Embed the tree (uniquely) in the plane so that
the clockwise order of the edge labels incident with every vertex is
increasing, and we complete the determination of $\phi^{-1}(T)$ straightforwardly
by inverting our planar dualization above.

\bigskip
	For example, for the tree $T_{1}$, with $n=11$, and edges
$47,37,23,39,13,15,56,5 \ 10,18,8 \ 11$, we find that

$$\phi^{-1}(T_{1})=((5 \ 6),(1 \ 6), (2 \ 3),(6 \ 9),(6 \ 7),(1 \
3),(9 \ 11),(3 \ 4),(7 \ 8),(9 \ 10)),$$

\bigskip
\noindent
as illustrated in Figure 3.

\begin{figure}[h]
\centering
\epsfig{file=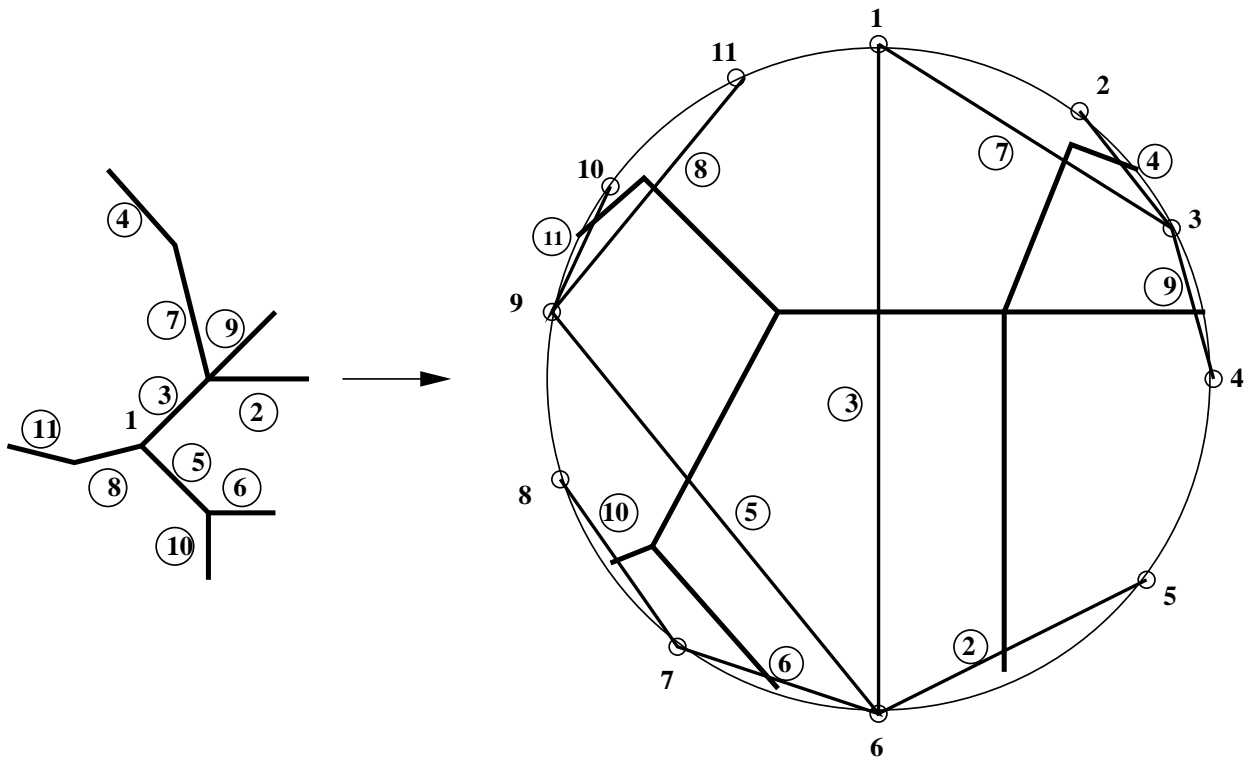,height=10cm}
\caption{Example of inverting the structural bijection.}
\label{figure_3}
\end{figure}

\bigskip
	As a remark, we mention that our circle chord diagram
construction leads directly to another bijection
between ${\mathcal F}_{n}$ and ${\mathcal T}_{n}$. Given a
circle chord diagram associated with
a factorization, simply ``push'' the edge labels in the unique
direction away from the vertex labelled 1. This gives an element of
${\mathcal T}_{n}$ and is clearly reversible. This bijection is
actually the same as that of Moszkowski~[\ref{M}], and the
description given above has appeared, independently, in
Poulalhon~[\ref{P}].

\section{Factorizations and Hurwitz numbers}

	The factorizations that we have considered in this paper
are special cases of a more general factorization problem in
 ${\mathcal S}_{n}$.  Consider $k$-tuples of
transpositions $(\sigma_{1},\ldots ,\sigma_{k})$ whose ordered product
 $\sigma_{1}\cdots\sigma_{k}$ is equal to an arbitrary permutation
 $\pi$, and such that the group generated
by $\sigma_{1},\ldots ,\sigma_{k}$ acts transitively on $[n]$. (When
 $\pi ={\mathcal C}_{n}$, as is the case in this paper, this transitivity
is forced. Note that, in general, transitivity means in
combinatorial terms that the graph on
vertices $[n]$, and edge $i \ j$ for each factor $(i \ j)$, is connected.)
For each $k$, the number of such factorizations is clearly constant on
the conjugacy class of $\pi$ in ${\mathcal S}_{n}$. Moreover, if
the conjugacy class has disjoint cycle distribution specified
by the partition $\alpha= (\alpha_1,\ldots ,\alpha_m )$ of $n$, with
 $m$ parts, then
the minimum choice of $k$ for which such factorizations exist
is $k=n+m-2$, and corresponding factorizations are called
{\em minimal transitive factorizations}. (Note that, from Goulden and
Jackson~[\ref{GJ}], there are exactly
 $m-1$ cuts in these factorizations, in addition
to $n-1$ joins.) The number of
these factorizations is given by
\begin{equation}\label{Hnum}
(n+m-2)!\,n^{m-3} \prod_{j=1}^m \frac{\alpha_j^{\alpha_j}}{(\alpha_j -1)!},
\end{equation}
from~[\ref{GJ}]. Such factorizations
arise in the study of ramified covers of the sphere by the sphere,
with branching
above infinity specified by $\alpha$, simple branching above other
specified points, and no other branching 
(see, for example~[\ref{A}],~[\ref{GJV}] and~[\ref{H}]).

\bigskip
	In the case where the number of factors is $k=n+m-2+2g$, for
an arbitrary nonnegative integer $g$, such factorizations
arise in ramified covers of the sphere by a surface of genus $g$. The
number of such covers, equal to the number of corresponding
factorizations as specified above, are called  Hurwitz numbers, and
are studied extensively in algebraic geometry (see, for example,
~[\ref{ELSV}] and~[\ref{FP}]). Here, the expression for $k$ is
a consequence of the Riemann-Hurwitz formula.

\bigskip
	Clearly the expression in~(\ref{Hnum}) specializes correctly to
give $n^{n-2}$ in the case that $\alpha =(n)$, for which $k=n-1$,
and these are the factorizations studied in this paper. We would
like to achieve a combinatorial understanding of this expression
for arbitrary $\alpha$. The bijection in this paper allows us
to specify transposition factors according to their difference index, and in
particular to identify those in which this difference equals $1$
(namely, the consecutive pairs). These pairs are mapped to leaves
in the corresponding tree, and the Pr\"ufer code bijection for
trees (see, e.g., Stanton and White~[\ref{SW}], p. 66) is based
on successively removing leaves of the tree, each
iteration yielding an element of $[n]$. Thus we can compose our
bijection with the Pr\"ufer bijection ``smoothly'', to identify
combinatorially each of the factors $n$ in the enumeration of
the factorizations of this paper. Now expression~(\ref{Hnum}) contains
many similar factors, and our hope is that the combinatorial
decompositions of this paper can be extended to explain these factors in the
general case. As a specific instance of this possible extension,
consider the case $\alpha =(n-1,1)$,
where expression~(\ref{Hnum}) becomes $(n-1)^n$, and the factorizations
would have a single cut. The simplicity
of this expression suggests that a nice combinatorial explanation of
the type referred to above should be possible, but we have not yet
been able to find one. 

\section*{Acknowledgements}
This work was supported by the Natural Sciences and Engineering
Research Council of Canada, through a grant to IG, and a PGSA
to AY. This work was partially completed while AY was visiting the
Field's Institute in Toronto. We thank Gilles Schaeffer for helpful
discussions.

\section*{References}

\begin{enumerate}
\item\label{A}
V.I. Arnol'd, {\em Topological classification of trigonometric polynomials
and combinatorics of graphs with an equal number of vertices and edges},
Functional Analysis and its Applications {\bf 30} (1996), 1--17.
\item\label{C}
A. Cayley, {\em A theorem on trees}, Quart. J., Oxford {\bf 23} (1889), 376--378.
\item\label{D}
J. D\'{e}nes, {\em The representation of a permutation as the product of
a minimal number of transpositions and its connection with the theory of
graphs}, Publ. Math. Institute Hung. Acad. Sci. {\bf 4} (1959), 63--70.
\item\label{ELSV}
T. Ekedahl, S. Lando, M. Shapiro and A. Vainshtein, {\em On Hurwitz
numbers and Hodge integrals}, C. R. Acad. Sci. Paris {\bf 328} S\'erie I
 (1999), 1171--1180.
\item\label{FP}
B. Fantechi and Pandharipande, {\em Stable maps and branch divisors},
math.AG/9905104
\item\label{GJ}
I.P. Goulden and D.M. Jackson, {\em Transitive factorizations into
transpositions and holomorphic mappings on the sphere}, Proc. Amer.
Math. Soc. {\bf 125} (1997), 51--60.
\item\label{GJV}
I.P. Goulden, D.M. Jackson and R. Vakil, {\em The Gromov-Witten potential
of a point, Hurwitz numbers and Hodge integrals}, Proc. London Math.
Soc. (to appear)
\item\label{GP}
I.P. Goulden and S. Pepper, {\em Labelled trees and factorizations of a
cycle into transpositions}, Discrete Math. {\bf 113} (1993), 263--268.
\item\label{H}
A. Hurwitz, {\em Ueber Riemann'sche Fl\"achen mit gegebenen Verzeigungspunkten},
Mathematische Annalen {\bf 39} (1891), 1--60.
\item\label{M}
P. Moszkowski, {\em A solution to a problem of D\'{e}nes: a bijection
between trees and factorizations of cyclic permutations}, European
J. Combinatorics {\bf 10} (1989), 13--16.
\item\label{P}
D. Poulalhon, ``Graphes et d\'{e}compositions de permutations'',
M\'{e}moire de DEA,
LIX, \'{E}cole Polytechnique, July 1997.
\item\label{SW}
D. Stanton and D. White, ``Constructive Combinatorics'', Springer-Verlag,
New York, 1986.
\end{enumerate}
\end{document}